\documentclass[12pt]{amsart}

\usepackage{amsmath,amsthm,amssymb,amscd,euscript,longtable,enumitem}

\usepackage{graphicx,color,mathtools,tikz}
\usepackage{hyperref,url}
\usepackage{nicematrix}

\usetikzlibrary{math,fit} 
\def \r{\mathbb R}
\def \rp{\mathbb{RP}}

\def \z{\mathbb Z}

\DeclareMathOperator{\GL}{GL}

\def \({\langle}
\def \){\rangle}

\usepackage[style=numeric-comp,natbib=true,backend=biber,maxnames=4,minnames=1]{biblatex}
\AtEveryBibitem{
	\clearfield{issn} 
	\clearfield{doi} 
	\clearfield{isbn} 
	
	\ifentrytype{online}{}{
		\clearfield{url}
	}
}

\renewbibmacro{in:}{}
\makeatletter
\def\keywords{\xdef\@thefnmark{}\@footnotetext}
\makeatother


\usepackage[todonotes={textsize=footnotesize}]{changes}
\definechangesauthor[color=blue]{OK}
\definechangesauthor[color=red]{MA}

\DeclareMathOperator{\sign}{sign}

\DeclareMathOperator{\sgn}{sgn}

\newtheorem{theorem}{Theorem}[section]

\newtheorem{proposition}[theorem]{Proposition}
\newtheorem{corollary}[theorem]{Corollary}
\theoremstyle{remark}
\newtheorem{remark}[theorem]{Remark}
\theoremstyle{definition}
\newtheorem{definition}[theorem]{Definition}
\newtheorem{example}[theorem]{Example}
\newtheorem{problem}{Problem}

\newtheorem{question}{Question}

\author[O.~Karpenkov]{Oleg~Karpenkov}
\address{Department of Mathematical Sciences\\ University of Liverpool\\ Peach Street \\ Liverpool L69~7ZL}
\email{karpenkov@liv.ac.uk}

\author[M.~Arnold]{Maxim~Arnold}
\address{Department of Mathematical Sciences\\ University of Texas at Dallas\\ 800, W. Campbell rd. \\ Richardson, TX, 75035}
\email{maxim.arnold@utdallas.edu}

\title{On Euclidean Algorithms for oriented linear Grassmanians}

\begin{date}  {\today} \end{date}
\begin{document}
	
\maketitle

\begin{abstract}
In this paper we study Euclidean algorithms and the corresponding continued fractions for oriented linear 
Grassmanians $G(k,n)$.
We propose two algorithms: Maximal
Element Elimination algorithm and 
Minimal Element Elimination algorithm.
The first algorithm reduces the absolute maximal value of the Pl\"ucker coordinates; the algorithm works only in $G(2,n)$. 
The second algorithm eliminates the 
Pl\"ucker coordinate with the smallest absolute values, while all other coordinates may increase; the algorithm works for arbitrary $G(2,n)$. We discuss basic features of these algorithms and formulate several natural open questions for further studies.
\end{abstract}


\section*{Introduction}

In this paper we introduce and discuss generalized subtractive Euclidean  
algorithms and the corresponding multidimensional continued-fraction 
sequences for the elements of the oriented integer linear Grassmanians $G(k,n)$ for $k\ge 2$.

\vspace{2mm}

To begin with, the case of Grassmanians $G(1,n)$ 
is very well studied since $G(1,n)$ is canonically identified with $\rp^n$.
For that reason, the subtractive algorithms on $G(1,n)$ coincide with the standard subtractive algorithms in Euclidean spaces.
One of the first subtractive algorithms for $\r^3$ is the classical Jacobi-Perron algorithm (see~\cite{Perron1907,Jacobi1868}), it was used to approach the problem of cubic periodicity by C.~Hermite in~\cite{Hermite1839}. 
Further subtractive algorithms were proposed by
V.~Brun~\cite{Brun1958} and
E.S.~Selmer~\cite{Selmer1961},
fully subtractive algorithm by F.~Schweiger~\cite{Schweiger1994}
and~\cite{Schweiger1995},
generalized subtractive algorithm~\cite{Schweiger1992},
 Tamura-Yasutomi algorithm~\cite{Tamura2009},
algebraic periodicity detecting algorithm~\cite{Karpenkov2022-H, Karpenkov2024-H},  etc.).
A nice collection of subtractive algorithms is in the book
by F.~Schweiger~\cite{Schw-2000}; see also the last chapter of~\cite{karpenkov-book-v2}. Geometric approach to multidimensional continued fraction theory were initiated 
by F.~Klein~\cite{Klein1895,Klein1896}, Minkowski~\cite{Minkowski1896} and Voronoi~\cite{voronoi1896}, see also 
in~\cite{Arnold2002,karpenkov-book-v2}.

\vspace{2mm}

The cases of Grassmanians $G(k,n)$ with $k\ge 2$ and $n-k\ge 2$ are more complicated, since such Grassmanians are not closed under standard addition operation on Pl\"ucker coordinates.  
Recall that Pl\"ucker coordinates
provide a natural embedding of the Grassmanian $G(k,n)$ into $\rp^N$ (where $N={{n}\choose{k}}$);
and that this embedding is not surjective for $k>1$. 
This develops a significant obstacle to straightforward generalizations of  algorithms to Grassmanians.
Indeed, the majority of the introduced Euclidean algorithms are based on vector subtraction,
while the difference of two vectors in $G(k,n)$ is not necessarily in $G(k,n)$ (with the only exceptions of $k=1$, and $k=n-1$).
As a consequence, both Euclidean subtractive and geometric algorithms have a limited usage for Grassmanians.

\vspace{2mm}

In this paper we develop several generalized subtractive Euclidean algorithms that preserve the structure of Pl\"ucker embeddings for Grassmanians. 
Our approach is to consider a vector reduction in $\r^n$ that evolves a plane in $G(k,n)$ to the first coordinate plane, i.e., to the plane for which only the first Pl\"ucker coordinate is non-zero.

\vspace{2mm}
In each of the algorithms we start with a collection of Pl\"ucker coordinates for some element of $G(k,n)$ (hence such coordinates satisfy the Pl\"ucker coordinate relations). We usually assume these coordinates to be non-zero. Our goal is to produce a sequence of subtraction-type reductions based on $GL(n,\z)$ transformations bringing  
these coordinates to a single non-zero integer coordinate.
The absolute value of this coordinate will be equal to the greatest common divisor of all the original Pl\"ucker coordinates.

\vspace{2mm}

Our first algorithm is the Maximal Element Elimination algorithm (see Subsection~\ref{Maximal Element Elimination}). It is limited to $G(2,n)$. It is based on the reduction of the Pl\"ucker coordinate with the maximal absolute value, which provide a uniform descent to the corresponding reduced form. 

\vspace{2mm}

The second algorithm is the Minimal Element Elimination algorithm (see Subsection~\ref{Minimal Element Elimination}), it works for all $G(k,n)$ . Here we reduce the Pl\"ucker coordinate with the minimal absolute value,
while the other coordinates may grow. Once one of the Pl\"ucker coordinates becomes zero we show how to reduce the dimension in order to continue in $G(k,n-1)$. 

\vspace{2mm}

The algorithms stop when they arrive to $G(k,k)$. The absolute value of only single coordinate here equals the greatest common divisor of the original Pl\"ucker coordinates.
Similarly to the classical Euclidean algorithm, the steps of the new algorithms can be considered as the elements of generalized continued fractions. Using them one arrives to a $k$-tuple of integer vectors, whose Pl\"ucker coordinates coincide with the original ones. 

\vspace{2mm}

\noindent
{\bf Organization of paper.}
In Section~\ref{Preliminary notions and definitions} we start with the classical definition of Pl\"ucker coordinates and discuss the problem of reduction to positive Pl\"ucker coordinates. We also recall the definitions of subtractive Euclidean algorithms
and study elementary subtractions in Grassmanians.
Further in Section~\ref
{Continued Fraction Algorithms for Grassmanians}
we introduce both Maximal Element Elimination algorithm for $G(2,n)$ and Minimal Element Elimination algorithm for $G(k,n)$ followed by discussion of their basic properties and an example.
We conclude the paper with a few directions for further studies in Section~\ref{Directions for further studies}.

\section{Preliminary notions and definitions}
\label{Preliminary notions and definitions}

A point or vector is said to be {\it integer} if all its coordinates are integers.
We say that a linear space is {\it integer} if it contains a full rank integer sublattice.

\vspace{1mm}

From now on we fix a basis of integer vectors and consider all vectors as a collection of coordinates in this basis.

\subsection{Pl\"ucker coordinates of Linear Grassmanians}

First of all denote by $G(d,n)$ the {\it linear Grassmanian} of $d$-planes in $\r^n$.

\vspace{2mm}

For a given $d$-dimensional space $\pi$ we consider 
its basis $(v_1,\ldots, v_n)$, and write it as a single $n{\times} d$-matrix: the rows of this matrix are the coordinates in the fixed basis of $\r^n$.
This matrix has precisely $N={n \choose d}$ minors of dimension $d$. The collection of such minors ordered 
lexicographically is said to be
the {\it Pl\"ucker coordinates} of $\pi$.
For a $d$-dimensional space $\pi\subset\r^n$ we denote them by
$P=(p_{i_1,\ldots, i_m})$, where all indices $(i_1,\ldots,i_m)$ satisfy $1\le i_1 < \dots <i_m\le n$.

\vspace{2mm}

The mapping of the elements of $G(d,n)$ to their Pl\"ucker coordinates is called  the
{\it Pl\"ucker embedding}
of $G(d,n)$ to $\rp^N$ with $N={n \choose d}$. Note that Pl\"ucker coordinates are uniquely defined up to a common scalar factor.

\subsection{Reduction to positive Pl\"ucker coordinates}
\label{Reduction to positive Plucker coordinates}

Let us start with the following natural definition.

\begin{definition}\label{Positive Skew Def}
We say that a $n$-tuple of linearly independent vectors in $\r^d$ are
{\it totally positive/non-negative}, if the corresponding Pl\"ucker coordinates 
are all positive/non-negative.
\end{definition}

Let us address the following question.

\vspace{2mm}

\noindent
\begin{question}
\label{que:positive} 
Does there exist an integer basis transformation sending a given collection of non-zero Pl\"ucker coordinates to positive ones?
\end{question}

\vspace{2mm}

The answer to this question is affirmative for $G(1,n)$ and for $G(2,n)$.
In the case of $G(1,n)$ one just swap the signs of the basis vectors corresponding to the negative Pl\"ucker coordinates. The case $k=2$ is a bit more complicated. 
We start with the following proposition.

\begin{proposition}\label{positivity-G2n}
Consider a linear space $\pi$ spanned by two vectors $v_1=(a_1,\ldots, a_n)$  and $v_2=(b_1,\ldots, b_n)$.
Then a pair $(v_1,v_2)$ is totally positive if and only if 
there exists a line $\ell\subset \r^2$ passing through the origin such that 
the 2-dimensional column vectors $(a_i,b_i)$ $i=1,\ldots, n$ are distinct
and are contained in one of the half-planes of $\r^2\setminus \ell$;
and they are ordered counterclockwise in this half-plane.
\end{proposition}

\begin{proof}
The Pl\"ucker coordinate $p_{i,j}$ equals the oriented volume of the parallelogram spanned by the vectors $(a_i,b_i)$ and $(a_j,b_j)$, which is positive by construction.
\end{proof}

\begin{remark}\label{swap}
Note that every $(2\times n)$ matrix can be transferred to a non-negative form by changing the  signs of vectors  $(a_i,b_i)$ for some $i$ and reordering these vectors afterwards.
Furthermore, one can send all the vectors to the half-plane $x_1\ge 0$.
\end{remark}

\begin{corollary}\label{ToPositiveInG(2,n)}
There exist an integer basis transformation sending a given col-
lection of non-zero Pl\"ucker coordinates in $G(2,n)$ to positive ones.
\end{corollary}

\begin{proof}
The proof is constructive, it is provided by the following algorithm.

\vspace{2mm}
\noindent
{\bf Input.} We are given by Pl\"ucker coordinates $P$. We assume that none of them 
are zeroes.

\vspace{2mm}
\begin{itemize}
\item {\bf Stage I.}
First of all we make all coordinates $p_{1,n}$ to be positive by changing the signs of the corresponding basis vectors in case when $p_{1,n}$ are negative.
Note that all the obtained vectors will be in one half-plane
with respect to the line passing through $(a_1,b_1)$ (bare in mind that we do not know the coordinates of the vector $(a_1,b_1)$, but only Pl\"ucker coordinates).

\vspace{2mm}
\item {\bf Stage II.}
It remains to swap the order of basis vectors according to the following rule: 
if $p_{i,j}<0$ then we swap $e_i$ and $e_j$.
We do this iteratively, taking always the first $p_{i,j}$ in the lexicographic order. 

In the worst case when all $p_{i,j}$ are negative 
(i.e. the vectors $(a_i,b_i)$ are all in the wrong order, cf.~Proposition~\ref{positivity-G2n}), we need to swap ${n-2 \choose 2}$ times.

\end{itemize}

\noindent{\bf Output.}
After the transformations of Stages~I and~II we have all vectors in one half-plane ordered counter clockwise.
According to Proposition~\ref{positivity-G2n} all Pl\"ucker coordinates will be positive.
\end{proof}

Note that the input of the above algorithm does not contain vectors $(a_i,b_i)$, but only the corresponding Pl\"ucker coordinates $P$.

\vspace{2mm}

The answer to the question \ref{que:positive} is affirmative for $G(3,5)$ as well since $G(3,5)$ is naturally isomorphic to $G(2,5)$.
However, for $G(3,6)$ the positivity cannot be achieved by reflections and reorderings only  due to the following counterexample.

\begin{example}\label{G(3,6)}
Recall that the elements of $G(3,6)$ has $20=\binom{6}{3}$ Pl\"ucker coordinates $p_{i_1,i_2,i_3}$; they correspond to all possible ordered triplets $\mathbf{i}=(i_1<i_2<i_3)$. Note that:

\begin{itemize}
\item Reflecting representative $v_j$ results in the sign change of all $p_\mathbf{i}$ with $j\in \mathbf{i}$. 
\item Swapping representatives $v_j$ and $v_k$ results in sign change of all $p_\mathbf{i}$ with $(j,k)\in \mathbf{i}$. 
\end{itemize}
Hence the former operation provides $10$ sign changes and the latter provides $4$ neither of which affect the parity of the amount of negative coordinates. 

\vspace{1mm}

\noindent
Consider the element of $G(3,6)$ defined by the following three vectors: 
$$\begin{bmatrix}
    1&0&0&1&1&1\\
    0&1&0&-3&-2&-1\\
    0&0&1&8&5&1
\end{bmatrix}
$$ 
Its Pl\"ucker coordinates are all positive except for $p_{4,5,6}$:
$$
P=(
1: 8: 5: 1: 3: 2: 1: 1: 5: 3: 1: 1: 1: 3: 7: 4: 1: 2: 1: -1
).
$$

Hence this element cannot be taken to totally positive by reordering or reflecting the vectors.  
\end{example}

\begin{remark}
Note that one can bring a given element of Grassmanian $G(d,n)$ to the positive form by the following iterative algorithm: 
 
\begin{enumerate}[label=(\roman*)]
    \item \label{step1} Fix first $(d-1)$ columns $v_1,\dots, v_{d-1}$ and reflect the columns $v_{m_1}$ with $m_1>(d-1)$ with $p_{1,\ldots,(d-1),m_1}<0$.

\item \label{step2} Next for any $m_1>m_2>(d-2)$ if $p_{1,\ldots,(d-2),m_2,m_1}<0$ the operation $v_{m_2}\mapsto v_{m_2}+v_{d-1}$ does not affect the coordinates from Step~\ref{step1}, increasing the coordinate $p_{1,\ldots,(d-2),m_2,m_1}$.

\item[{(k)}] \label{stepk} Next for any $m_1>m_2>\dots>m_k>(d-k-1)$ if $p_{1,\ldots,(d-k-1),m_k,\ldots, m_1}<0$ the operation $v_{m_k}\mapsto v_{m_k}+v_{d-k-1}$ does not affect the coordinates from the previous steps, increasing the coordinate $p_{1,\ldots,m_k,\dots,m_1}$.
\end{enumerate}
In this method we used a rather arbitrary $SL(n,\z)$ transformations.
Despite the case of symmetries about coordinate planes and swapping coordinates, the new ones can change the maximal absolute value of Pl\"ucker coordinates. By that reason such method cannot be used for Maximal Element Elimination algorithm.
Thus the question which remains open to our knowledge can be stated as
\begin{problem}\label{prb:positivity}
     Is there some natural way to make all Pl\"ucker coordinates positive without changing their maximal absolute value?
\end{problem}
\end{remark}

\subsection{Integer linear Grassmanians}
Let $\z G(d,n) \subset G(d,n)$ denote the subset of all integer linear spaces. Note that $\z G(d,n)$ is everywhere dense in $G(d,n)$ in the standard topology.

\begin{proposition}
A linear space $\pi \in G(d,n)$ is integer if and only if its Pl\"ucker coordinates define a rational point in $\rp^{N-1}$.
\end{proposition}

Recall that the {\it integer volume} of the sublattice generated by $v_1,\ldots, v_d\in \pi$ is the index of this sublattice in the integer lattice of $\pi$.

\begin{remark}
Consider linearly independent integer vectors $v_1,\ldots, v_d$ spanning $\pi$. 
\\
Then the greatest common divisor of the Pl\"ucker coordinates for these vectors coincides with the integer volume of the sublattice generated by $v_1,\ldots, v_d$,
see e.g., in~\cite{karpenkov-book-v2}.
\end{remark}

We say that a vector of Pl\"ucker coordinates is {\it primitive integer}
if all its entries are integers and their greater common divisor equals $1$.

Integer Pl\"ucker coordinates are defined up to multiplication by $\pm 1$ of all elements simultaneously.

\begin{remark}
Let us say a few words about the simplest case:
the case of lines through the origin in the plane. Let a line $\ell$ contain some vector $(a,b)$. Then its Pl\"ucker coordinates is a point  $(a:b)\in \rp^1$. Here the ratio $b/a$ determine the tangent of the slope of $\ell$. In case of an integer line the tangent of the slope $b/a$ is rational. 
\end{remark}

\subsection{Classic subtractive algorithms in $\r^n$}

Let us say a few words regarding classic subtractive algorithms in $\r^n$. They will be substantially used in the below algorithms for Grassmanians.

\subsubsection{Euclidean algorithm}

We start with the Euclidean algorithm. Given two integers $(p_0,q_0)$ we perform the following operation:
$$
(p_i,q_i) \to (p_{i+1},q_{i+1})=
(q_i, p_i-a_iq_i),
\quad 
\hbox{where $a_i=\lfloor p_i-q_i\rfloor$}.
$$
After a finitely many iterations (let it be $N+1$), the algorithm will arrive to the pair of integers $(\gcd(p_0,q_0),0)$ where it terminates.

\vspace{2mm}

Let us collect the output $(a_0,\ldots, a_N)$ of the above algorithm. 
This output forms a {\it regular continued fraction} of $p_0/q_0$. Namely
$$
\frac{p_0}{q_0}
=a_0+\frac{1}{\displaystyle a_1+\frac{1}{\displaystyle
a_2+\frac{1}{\displaystyle\ddots+\frac{\displaystyle 1}{a_N}}}}.
$$
It is usually denoted by $[a_0;a_1:\ldots:a_{N}]$.

Note that the Euclidean algorithm can be represented in terms of $\GL(2,\z)$ matrix multiplication form as follows:
$$
\left(
\begin{matrix}
\gcd(p_0,q_0)\\
0\\
\end{matrix}
\right)
=
\left(\prod\limits_{i=1}^N
\left(
\begin{matrix}
0&1\\
1&-a_i \\
\end{matrix}
\right)
\right)
\left(
\begin{matrix}
p_0\\
q_0\\
\end{matrix}
\right)
$$
The product of matrices on the right hand side correspond to elementary lattice basis transformation reducing the coordinates of the starting vector to $\gcd(p_0,q_0)$.  

\subsubsection{Multidimensional generalisations}

Assume that we have constructed a vector $(x_{1,i},\ldots, x_{n,i})$ and let $x_{2,i}\ne 0$.
Then the next vector $(x_{1,i+1},\ldots, x_{n,i+1})$ is set as follows
$$
\Big(
x_{2,i},
x_{3,i}-\Big\lfloor\frac{x_{3,i}}{x_{2,i}}\Big\rfloor {x_{2,i}}, 
\ldots,
x_{n,i}-\Big\lfloor\frac{x_{n,i}}{x_{2,i}}\Big\rfloor {x_{2,i}},
x_{1,i}-\Big\lfloor\frac{x_{1,i}}{x_{2,i}}\Big\rfloor {x_{2,i}}\Big).
$$
As an output of this step we have the following pair of integers:
$$
\Big(
\Big\lfloor\frac{x_{3,i}}{x_{2,i}}\Big\rfloor,
\ldots,
\Big\lfloor\frac{x_{n,i}}{x_{2,i}}\Big\rfloor,
\Big\lfloor\frac{x_{1,i}}{x_{2,i}}\Big\rfloor
\Big).
$$
It is called the {\it $i$-th element} of the Jacobi-Perron multidimensional continued fraction algorithm.

\vspace{2mm}

Once $x_{2,N}=0$, one swaps the second and the last coordinates and continue the process with one coordinate less. The algorithm terminates when there is only one non-zero coordinate left. 

\begin{remark}
Note that the $i$-th step of the Jacobi-Perron multidimensional continued fractions algorithms are described by multiplication of vectors by the following $\GL(n,\z)$ matrix:
$$
\left(
\begin{matrix}
0&1&0&0& \ldots &0\\
0&-\frac{x_{3,i}}{x_{2,i}}&1&0& \ldots &0 \\
0&-\frac{x_{4,i}}{x_{2,i}}&0&1& \ldots &0 \\
\vdots&\vdots& \vdots& \vdots& \ddots & \vdots\\
0&-\frac{x_{n,i}}{x_{2,i}}&0&0& \ldots &1 \\
1&-\frac{x_{1,i}}{x_{2,i}}&0&0& \ldots &0 \\
\end{matrix}
\right).
$$
\end{remark}

\begin{remark}
In the multidimensional case ($n\ge 3$) there is no unique choice of the subtractive algorithms.
In this paper we consider the subtractive algorithm 
that in the case $n=3$ coincides with one of the oldest subtractive algorithms, the Jacobi-Perron algorithm 
(see~\cite{Perron1907,Jacobi1868}).
Instead of Jacobi-Perron algorithm one can use any other subtractive algorithms (see e.g., in~\cite{Brun1958,Selmer1961,Schweiger1994,Schweiger1995,Schweiger1992,Tamura2009,Karpenkov2022-H, Karpenkov2024-H},  etc.).
\end{remark}

\subsection{Elementary subtraction of basis vectors}
\label{Elementary subtraction of basis}

\begin{definition}
{\bf(Elementary subtraction of basis vectors.)}
Let $1\le s\ne t \le n$.
The linear map defined on basis vectors of $\r^n$ as  
$$
e_i \to 
\left\{
\begin{array}{ll}
e_i, &\hbox{if $i\ne s$};\\
e_s-e_t,& \hbox{otherwise}\\
\end{array}
\right.
$$
is said to be an {\it elementary reduction}.
Denote it by $T_{s,t}$.
\end{definition}

\begin{definition}
Let $\mathbf{i}$ be an unordered $k$-tuple of distinct positive integers. Denote by $\sigma(\mathbf{i})$ the permutation that puts the indices to the increasing order.

\vspace{2mm}

\noindent
{
Let $s,t$ be two positive integers such that $s \in \mathbf{i}$ and $t \in \mathbf{i}$. By $\mathbf{i}_{s\to t}$ we denote the $k$-tuple of integers with $s$ replaced by $t$
(without changing the order).    
}
\end{definition}

\begin{proposition}
Consider an ordered $k$-tuple $\mathbf{i}$.
The subtraction $T_{s,t}$ changes the coordinates $p_\mathbf{i}$ as follows
$$
p_\mathbf{i} \to 
\left\{
\begin{array}{ll}
p_{\mathbf{i}}-
\sign(\sigma(\mathbf{i}_{s\to t}))\cdot
p_{\sigma(\mathbf{i}_{s\to t})},& \hbox{if $s\in \mathbf{i}$ and $t\notin \mathbf{i}$};\\
p_{\mathbf{i}}, &\hbox{otherwise}.\\
\end{array}
\right.
$$
\end{proposition}

\begin{proof}
The proof follows from the direct formula for determinants.
\end{proof}

For the case $G(2,n)$ the above expression for $T_{s,t}$ takes the form
\begin{equation}
    \label{eq:T_st}
     \begin{cases} p_{j,s}\mapsto
     p_{j,s}+\sgn(t-j)p_{j,t}\\
     p_{s,j}\mapsto p_{s,j} +\sgn(t-j) p_{t,j}
     \end{cases}\quad \mbox{ for } j\ne t
\end{equation} 
and thus can be regarded as column operation on the matrix of Pl\"ucker coordinates, tweaked in a way to preserve skew-symmetry.

\tikzset{highlight/.style={rectangle,
			fill=red!15,
			blend mode = multiply,
			rounded corners = 0.5 mm,
			inner sep=1pt}}
		\tikzset{lowlight/.style={rectangle,
				fill=blue!15,
				blend mode = multiply,
				rounded corners = 0.5 mm,
				inner sep=1pt}}

\[P=\begin{pNiceMatrix}[code-after = {\tikz { \node[highlight, fit = (4-7) (6-7)] {};
\node[highlight, fit = (4-1) (4-9)] {};
\node[highlight, fit = (1-4) (9-4)] {};
\node[lowlight, fit = (1-7) (9-7)] {};}}]
0&p_{1,2}&\Cdots&p_{1,s}&\Cdots&&p_{1,t}&\Cdots&p_{1,n}\\
\Vdots&\Ddots&\Ddots&\Vdots&&\longleftarrow&\Vdots&&\Vdots\\
&&0&p_{s-1,s}&&&\\
-p_{1,s}&\Cdots&&0&p_{s,s+1}&\Cdots&p_{s,t}&\Cdots&p_{s,n}\\
\Vdots&&&\Vdots&\Ddots&\Ddots&\Vdots&&\Vdots\\
-p_{1,t}&\Cdots&&&\Cdots&0&p_{t,t+1}&\Cdots&p_{t,n}\\
\Vdots&&&&&&0&\Ddots&\Vdots\\
&&&&&\longleftarrow&\Vdots&\Ddots&p_{n-1,n}\\
-p_{1,n}&\Cdots&&-p_{s,n}&\Cdots&&-p_{t,n}&\Cdots&0\\
\end{pNiceMatrix}\]

\begin{remark}
The subtractions $T_{s,t}$ and $T_{t,s}$ lead to the same set of absolute values of the coefficients, in particular their maximum is the same.   
\end{remark}

\begin{remark}
We always can make our maximal element to be not at the second diagonal $(j,j+1)$.
Namely, if $p_{i,i+1}$ is the maximal element
for the positive Grassmanian $\{v_1,v_2,\ldots,v_n\}$, then $p_{1,n}$ is the maximal element for the positive Grassmanian $\{v_{i+1},v_{i+2},\ldots,v_n,-v_1,-v_2,\ldots,-v_i\}$. In such settings we always can use subtraction $T_{j,j-1}$:

\end{remark}

\section{Continued Fraction Algorithms for Grassmanians}
\label{Continued Fraction Algorithms for Grassmanians}

Consider a $k$-tuple of vectors in $\r^n$: $(v_1,\ldots, v_k)$. 
In what follows we denote by $w_i\in\r^k$ the vectors of $i$-th coordinates of $v_1,\ldots, v_k$ for $i=1,\ldots, n$. 
Note that the matrices of $v$-vectors and $w$-vectors are transpose to each other.   

\vspace{2mm}

In what follows we will use the lexicographic order for the set of all indices of Pl\"ucker coordinates.

\subsection{Input and Output of the algorithms}
In this subsection we discuss the common items for both algorithms: input data and their outputs.

\vspace{2mm}

\noindent{
{\bf Input:}} Let $k<n$ be two positive integers.
We are also given by a collection $P$ of $N={{n}\choose{k}}$ integers, that defines an element of $G(k,n)$.

\vspace{2mm}

\noindent{
{\bf Tasks of the algorithm:}}
\begin{itemize}
    \item The first task of the algorithm is to find $k$-tuple of linearly independent integer vectors, whose Pl\"ucker coordinates coincide with $P$.
\vspace{1mm}
    
    \item The second task is to build a ``generalized continued fraction'' that encodes the basis transformations of the algorithm leading to $k$-tuple of vectors of the first task.  
\end{itemize}

\begin{remark}
Both tasks of the proposed algorithms are in the spirit of the classic Euclidean algorithm that 
finds the $\gcd$ of a pair of integers.
New algorithms provide a $k$-tuple of integer vectors together with the index of sublattice they generate in the $k$-plane they span.
This sublattice index can be seen as a straightforward generalization of $\gcd$.

\vspace{2mm}

\noindent{
The outcome of the second task is the sequence of $\GL(n,\z)$ transformations, similar to the ones generated by classical subtractive multidimensional continued fraction algorithms (see e.g., in~\cite{Schw-2000} and Section~27.4 of~\cite{karpenkov-book-v2}).}
\end{remark}

\vspace{2mm}
\noindent{
{\bf Output:}}
As an output we have the following:

\begin{itemize}
    \item A {\bf $k$-tuple of integer vectors} whose Pl\"ucker coordinates equal $P$. (In addition we will get the index of sublattice they generate.) 

    \item A {\bf sequence of $\GL(n,\z)$ transformations}, that we treat as the {\it continued-fraction sequence} for the algorithm. 
\end{itemize}

As in the classic Euclidean algorithm, the resulting vectors (together with the index of sublattice they generate) are obtained at the end of the algorithm, while the elements of the continued-fraction sequence are generated at the steps of the algorithm.

\begin{remark}
Both algorithms, applied to non-integer input $P$, provide the Grassmanian analog of the classical Gauss map. 
We encourage to study the properties of the corresponding dynamic systems.
\end{remark}

\subsection{Maximal Element Elimination algorithm for $G(2,n)$}
\label{Maximal Element Elimination}
Let us describe an algorithm that works for the case of $G(2,n)$.

\subsubsection{Annulation Step}
The purpose of this step is annulate one of the Pl\"ucker coordinates.
In the case of maximal element algorithm we iteratively reduce the coordinate with the maximal absolute value using integer bases transformations in $\r^n$ until we hit a zero Pl\"ucker coordinate.

\vspace{2mm}

\begin{itemize}
\item
{\noindent 
{\bf Stage I.}} Produce a basis where all the Pl\"ucker coordinates $p_{i,j}$ (with $i<j$) are positive. 
This is done by the algorithm in the proof of Corollary~\ref{ToPositiveInG(2,n)}.
The corresponding integer transformation matrix is
{\it added to the continued-fraction sequence}.

\vspace{2mm}

\item
{\noindent 
{\bf Stage II.}}
Select a Pl\"ucker coordinate $p_{i,j}$ (with $i<j$) attaining the maximal value.
In case if there are several such elements, we pick the one with the smallest lexicographic order of the index with one exception: we always prioritize $p_{1,n}$.

\vspace{2mm}

\item
{\noindent 
{\bf Stage III.}}
This stage is performed only if in Stage~II the coordinate $p_{i,i+1}$ was selected, i.e., the case $j=i{+}1$.
First, we cyclically shift the coordinate vectors in $\r^n$ by $n-i$ positions. Secondly, we multiply all the basis vectors of $\r^n$ from $i$ to $n$ by $-1$.
As a result, all the Pl\"ucker coordinates $p_{i,j}$ (with $i<j$) are positive and 
the coordinate $p_{1,n}$ becomes one of the maxima.
We select it for the next stage.

The corresponding basis transformation matrix is 
{\it added to the continued-fraction sequence}.

\vspace{2mm}

\item
{\noindent 
{\bf Stage IV.}}
Here we get the maximal coordinate $p_{i,j}$ with $i<j-1$. Now we perform the subtraction
$$
e_j\to e_j-e_{j-1}.
$$
The corresponding basis transformation matrix is 
{\it added to the continued-fraction sequence}.
\\
Until all Pl\"ucker coordinates are non-zero, then we {\bf return to Stage~I}. 

\vspace{2mm}

\item{
{\bf Go to Dimension Reduction Step.}
Once we get some zero Pl\"ucker coordinate $p_{s,t}=0$ for some $s<t$ we terminate the iterations and move to Dimension Reduction Step below. 
}
\end{itemize}

\vspace{2mm}

\subsubsection{Dimension Reduction Step}

After the annulation step we are in the situation of some Pl\"ucker coordinate $p_{s,t}$ being equal zero.

\vspace{2mm}

\begin{itemize}
\item
{\noindent 
{\bf Stage I.}} Let us swap the coordinates of $\r^n$ such that in new coordinates $p_{1,2}=0$.

\item
{\noindent 
{\bf Stage II.}}
We have two cases here.

\vspace{1mm}

\subitem{(1)} 
Assume that there is some $i$ such that $p_{1,i}\ne 0$. Then the coefficient of proportionality of the vectors $w_1$ and $w_2$ is precisely $p_{2,i}:p_{1,i}$.
We apply the usual two-dimensional subtractive Euclidean algorithm to the pair $(v_1,v_2)$ in order to annulate $v_1$. Here we make the  subtractions following the Euclidean algorithm for the pair $(p_{1,i}, p_{2,i})$ until we get the $p_{1,i}$ to be a zero.
Once this is done, the plane under consideration is in the plane $x_1=0$ in the current basis. 
At this point we exclude the first coordinate. This reduces the dimension by one.

\vspace{1mm}

\subitem{(2)}
Assume that that for all $i$ we have $p_{1,i}= 0$. This implies that $w_1=0$. So we exclude the first coordinate.

\vspace{2mm}

\item{
{\bf Go to Annulation Step or Terminate.}
Once we eliminated one of the coordinates, we move to the Annulation Step with one dimension less. 
In case if we reduced the dimension to $k$, we terminate. 
}
\end{itemize}

\vspace{2mm}

The corresponding basis transformation matrices of these stages are 
{\it added to the continued-fraction sequence}.

\begin{remark}\label{slnz-remark-1}
For simplicity we work with $\GL(n,\z)$: we send the excluded coordinate to the last position.
The further algorithm will work in dimensions less then $n$, so all the lower-dimensional transformations will be supplemented by the identity transformation for the excluded coordinates.
\end{remark}

\subsection{Minimal Element Elimination algorithm for $G(k,n)$}
\label{Minimal Element Elimination}
The above algorithm does not work 

\subsubsection{Annulation step}
On the contrary to the Maximal element elimination algorithm, in this algorithm we do not intend to make all the Pl\"ucker coordinates positive.

\vspace{2mm}

\begin{itemize}
\item{
{\bf Stage I.}
Select a Pl\"ucker coordinate $p_{i_1,\ldots, i_k}$ (for $i_1<\ldots<i_k$) with the minimal  absolute value. In case of several minimal absolute values we take the one with the smallest lexicographical index.
}

\vspace{2mm}

\item{
{\bf Stage II.}
Select a Pl\"ucker coordinate $p_{j_1,\ldots, j_k}$ (for $j_1<\ldots<j_k$) that $j_s=i_s$, except for $s=t$.
The values of $t$ and $j_t$ are defined as follows. 

\vspace{1mm}

\noindent
$(${\it i}$)$ $t$ is the minimal index such that whether the index $i_{t}-1$ (considered only if $t>1$) or the index $i_{t}+1$ (considered only if $t<n$) are not in the set $\{i_1,\ldots, i_k\}$;

\vspace{1mm}

\noindent
$(${\it ii}$)$ We select $j_t=i_{t}{-}1$ if this option is available and $j_t=i_{t}{+}1$ otherwise.
}

\vspace{2mm}

\item{
{\bf Stage III.} Subtract
$$
e_{i_t} \to e_{i_t}-\bigg\lfloor\frac{p_{j_1,\ldots, j_k}}{p_{i_1,\ldots, i_k}}\bigg\rfloor e_{j_t}.
$$
The corresponding basis transformation matrix is 
{\it added to the continued-fraction sequence}.
}
\\
If all Pl\"ucker coordinates are non-zero, then we {\bf return to Stage~I}. 

\vspace{2mm}

\item{
{\bf Go to Dimension Reduction Step.}
Once we get some zero Pl\"ucker coordinate $p_{s,t}=0$ for some $s<t$ we terminate the iterations and move to the dimension reduction step. 
}
\end{itemize}

\subsubsection{Dimension reduction step}
This is the multidimensional analog of the Dimension reduction for $G(2,n)$ described in the Maximal element algorithm.

\vspace{2mm}

\begin{itemize}
\item
{\noindent 
{\bf Stage I.}} Swap the coordinates of $\r^n$ such that in new coordinates 
$$
p_{1,\ldots, k}=0.
$$

\vspace{2mm}

\item
{\noindent 
{\bf Stage II.}}
Let us find a nonzero Pl\"ucker coordinate $p_{1,\ldots, s,i_{s+1},\ldots, i_{k}}$ such that

--- all the coordinates $p_{1,\ldots, s, s+1, j_{s+2},\ldots, j_{k}}$ are zero;

--- $(i_{s+1},\ldots, i_{k})$ is lexicographically the first $(k-s)$-tuple that provides a non-zero Pl\"ucker coordinate.
\\
Note that $s<k$.

\vspace{2mm}

\item
{\noindent 
{\bf Stage III.}}
From the first condition of Stage II we get that the vectors $v_1,\ldots, v_{s+1}$ are linearly dependent.
Moreover they satisfy the following linear equation with integer coefficients:
$$
\left(\sum\limits_{j=1}^s
p_{1,\ldots,j-1,s+1,j+1,\ldots,s,i_{s+1},\ldots, i_k}
w_j
\right)
-
p_{1,\ldots,s,i_{s+1},\ldots, i_k}w_{s+1}
=0.
$$
Denote their coefficients by $(a_1,\ldots, a_{s+1})$.

\vspace{2mm}

\item{
{\bf Stage IV.}
Here we again have two cases. 

\vspace{1mm}

\subitem{(1)} Let $(a_1,\ldots, a_{s+1})=(0,\ldots, 0)$.
Then $w_{s+1}$ is a zero vector.
We swap $e_1$ and $e_{s+1}$ to get $w_1=0$.

\vspace{1mm}

\subitem{(2)} Let $a_i\ne 0$ for some $1\le i\le s+1$.

Then we apply the Jacobi-Perron algorithm to $(a_1,\ldots ,a_{s+1})$ that will bring the vector to $(a,0,\ldots,0)$.
As a result, in the new basis we have $aw_1=0$, and hence $w_1=0$.

\vspace{2mm}

\noindent
In both cases we exclude the first coordinate.
}

\vspace{2mm}

\item{
{\bf Go to Annulation Step or Terminate.}
Once we eliminated one of the coordinates, we move to the Annulation step with one dimension less. 
In case if we reduced the dimension to $k$, we terminate. 
}
\end{itemize}

The corresponding basis transformation matrices of these stages are 
{\it added to the continued-fraction sequence}.

\begin{remark}\label{slnz-remark-2}
Similar the case of Maximal element elimination algorithm (see Remark~\ref{slnz-remark-1}) 
we intend to work with $\GL(n,\z)$ basis transformations.
By that reason we send the excluded coordinate to the last position.
The further algorithm will work in dimensions less then $n$, so all the lower-dimensional transformations will be supplemented by the identity transformation for the excluded coordinates.
\end{remark}

\begin{remark}
Jacobi-Perron algorithm can be replaced to any other subtractive algorithms. For subtractive algorithms we refer to~\cite{Schw-2000} and Section~27.4 of~\cite{karpenkov-book-v2}). 
\end{remark}

\subsection{Output of the algorithms}

{\bf Termination of the algorithm.}
The algorithms performs $n-k$ dimension reductions ($k=2$ in the case of Maximal element elimination algorithm ). Once we arrive to $G(k,k)$ we terminate.
At this stage we have 
\begin{itemize}
\item $\hat p=\hat p_{1,\ldots,k}$ the only Pl\"ucker coordinate of $G(k,k)$;

\item the continued-fraction sequence $(A_i)_{i=1}^N$ of $\GL(n,\z)$ transformations.
\end{itemize}

Finally to complete the output we should produce vectors whose Pl\"ucker coordinates are as in the input.

\begin{itemize}
    \item The basis of the sublattice $(\hat v_1,\ldots,\hat v_k)$ with Pl\"ucker coordinates $P$ are obtained from the basis vectors $\hat p e_1, e_2\ldots, e_k$ by inverse multiplication on all matrices of the continued-fraction sequence:
\begin{equation}\label{eq2}
\hat v_1=\Big(\prod\limits_{i=1}^n (A_{i})^{-1}\Big)\hat p e_1, \qquad  
   \hat v_i=\Big(\prod\limits_{i=1}^n (A_{i})^{-1}\Big)e_i
   \quad \hbox{for $i=2,\ldots, k$.}
\end{equation}    
    Here when we increase the dimension, we simply add zero coordinates to all vectors. 
\end{itemize}

\vspace{2mm}

\noindent
{\bf Output:} $\Big((\hat v_1,\ldots,\hat v_k),(A_i)_{i=1}^N\Big)$.

%


\vspace{2mm}

\noindent
{\bf Finiteness termination time:} 
Each Dimension Reduction Step and and each Annulation Step of the algorithm are performed in a finitely many operations. It is sufficient to do $n-k$ iterations of pairs of these steps in order to arrive to the termination state.

\begin{remark}
If the input is irrational then the running time of the algorithm can be infinite, in this case we will have an infinite continued fraction as an output.
It is essential to study periodic algorithms (Lagrange's theorem), distribution of partial quotients (Gauss-Kuzmin statistics), and further properties of the corresponding generalized Gauss map.
\end{remark}

\subsection{Basic properties of the algorithms}

Finally we have the following property.

\begin{theorem}\label{plucker-volume}
The Pl\"ucker coordinate $\hat p$ is the volume of the sublattice generated by the vectors $(\hat v_1,\ldots,\hat v_k)$. 
\end{theorem}

\begin{proof}
On the one hand, the volume of the sublattice is invariant under all the basis transformation of the step of the algorithm.
On the other hand the index of the $k$-dimensional sublattice in $\mathbb Z^k$ coincides with the absolute value of its Pl\"ucker coordinate (i.e., the determinant
of the basis for the sublattice). 
\end{proof}

\begin{theorem}\label{plucker-volume}
For a give set of integer Pl\"ucker coordinates $P$ both algorithms generate a $k$-tuple of lattice vectors that coincide with $P$.
\end{theorem}

\begin{proof}
This follows from the construction. Here we have a sequence of integer linear transformations (see Equation~\ref{eq2}) bringing 
$\hat p e_1, e_2\ldots, e_k$ to the required basis.
\end{proof}

Note that the Maximal Element Elimination algorithm works only for the case of $G(2,n)$ due to the discussion in the section \ref{Reduction to positive Plucker coordinates} and Example~\ref{G(3,6)} above.
So we arrive to the following reformulation of the Problem \ref{prb:positivity} 

\begin{problem}
Generalize the Maximal Element Elimination algorithm to the case of $G(k,n)$ for $k\ge 3$.
\end{problem}

\subsection{Admissible $k$-tuple of vectorts}

{\it Which integer vectors have the same Pl\"ucker coordinates?} 
The answer to this question is rather simple. 
Let $\hat p$ be the last Pl\"ucker coordinate obtained as the first output of the algorithm. 
These are the $k$-tuples of integer linearly independent vectors that generate sublattice on the integer lattice in the $k$-plane of index $\hat p$. The Pl\"ucker coordinates of these $k$-tuples will be as desired (up to simultaneous choice of the sign).   

\vspace{2mm}

Let us indicate how to construct them.
In the output of the algorithms we have chosen the vectors $(\hat p e_1, e_2\ldots, e_k)$ and applied to them the elements of the constructed continued-fraction sequence.
Instead of taking $(\hat p e_1, e_2\ldots, e_k)$ one can choose any other integer basis of $\r^k$ with determinant $\hat p$.
Equivalently one might say that we start with an integer matrix of determinant $\hat p$.
Further we again apply the elements of the constructed continued-fraction sequence (as in Equation~\ref{eq2}). This will generate the set of integer $k$-tuples in $\r^n$ with the given Pl\"ucker coordinates.

\subsection{Working example}
\label{Subtractive Euclidean algorithm in Grassmanians}

Recall that there is a natural bijection between the set of all rank-2 skew-symmetric matrices in $PSL(n,R)$ and the Grassmannian $G(2,n)$. 
For instance, thanks to the Pl\"ucker relations
$$p_{i,k}p_{j,l} = p_{i,j}p_{k,l} + p_{jk}p_{ki},\qquad  \mbox{ for } i < j < k < l$$
one can consider any $G(2,n)$ as decorated Frieze pattern. Thus, for example if $p_{j,j+1}=1$ for all $j$ , one recovers usual Frieze (see~\cite{FREEZE} for more on friezes and \cite{CUNTZFrieze} for decorated friezes.) In this regard, our algorithms correspond to the natural operations on the corresponding cluster variety. Here we illustrate our maximal element elimination algorithm for the first non-trivial case of $G(2,4)$.

\begin{example}
Consider the following set of Pl\"ucker coordinates:
$$(p_{1,2},p_{1,3},p_{1,4},p_{2,3},p_{2,4},p_{3,4})=
(10,10,12,-15,3,21)
$$

Then after the reflections and reordering step, we get 
$(10,10,12,15,21,3)$.
Maximal Element Elimination algorithm then proceeds as follows:
$$\begin{pmatrix}
    0&10&10&12\\
    -10&0&15&21\\
    -10&-15&0&3\\
    -12&-21&-3&0
\end{pmatrix}\stackrel{T_{3,4}}{\mapsto} \begin{pmatrix}
    0&10&10&2\\
    -10&0&15&6\\
    -10&-15&0&3\\
    -2&-6&-3&0
\end{pmatrix}$$

Maximal element is now at the position $(2,3)$. After cyclic rotation  $(v_1,v_2,v_3,v_4)\mapsto (-v_3,-v_4,v_1,v_2)$ the maximal element is now at the position $(1,4)$ and so the subtraction $T_{3,4}$ reduce it.

$$\begin{pmatrix}
    0&3&10&15\\
   -3&0&2&6\\
    -10&-2&0&10\\
    -15&-6&-10&0
\end{pmatrix}\stackrel{T_{3,4}}{\mapsto}
\begin{pmatrix}
    0&3&10&5\\
   -3&0&2&4\\
    -10&-2&0&10\\
    -5&-4&-10&0
\end{pmatrix}$$

Now the maximal elements are at the positions $(1,3)$ and $(3,4)$. We use the \emph{lexicographic order for the indices modulo (n)}, i.e. $(i_1,j_1)\prec (i_2,j_2)$ if $(i_1\mod(n)<i_2\mod(n))$ or $i_1=i_2$ and $j_1\mod(n)<j_2\mod(n)$ , hence the first maximum is at the position $(1,3)$.
Thus applying $T_{2,3}$ we get

$$\begin{pmatrix}
    0&3&7&5\\
   -3&0&2&4\\
    -7&-2&0&6\\
    -5&-4&-6&0
\end{pmatrix}
\stackrel{T_{2,3}T_{3,4}T_{2,3}}{\mapsto}
\begin{pmatrix}
    0&3&1&1\\
   -3&0&2&2\\
    -1&-2&0&0\\
    -1&-2&0&0
\end{pmatrix}.$$

This is the end of the first stage of algorithm, as $v_4=v_3$. Hence, we could project on three first coordinates and continue.
$$
\begin{pmatrix}
    0&3&1\\
   -3&0&2\\
    -1&-2&0
\end{pmatrix}.\stackrel{(1,2,3)\to (-2,-3,1)}{\mapsto}
\begin{pmatrix}
    0&2&3\\
   -2&0&1\\
    -3&-1&0
\end{pmatrix}\stackrel{T_{2,3}}{\mapsto}
\begin{pmatrix}
    0&2&1\\
   -2&0&1\\
    -1&-1&0
\end{pmatrix}
$$

$$\stackrel{(1,2,3)\to (-2,-3,1)}{\mapsto}\begin{pmatrix}
    0&1&2\\
   -1&0&1\\
    -2&-1&0
\end{pmatrix}\stackrel{T_{2,3}T_{2,3}}{\mapsto} \begin{pmatrix}
    0&1&0\\
   -1&0&1\\
    0&-1&0
\end{pmatrix}\stackrel{(1,2,3)\to (-2,-3,1)}{\mapsto}\begin{pmatrix}
    0&1&1\\
   -1&0&0\\
    -1&0&0
\end{pmatrix}
$$

We arrived to a single coordinate $(p)=(1)$ in $G(2,2)$. Now applying the inverse transformations to the pair of vectors $(p\cdot e_1,e_2)=(e_1,e_2)$ we arrive to a basis $(v_1,v_2)$ of the integer plane in $G(2,4)$:
$$
\begin{pmatrix}
    v_1\\v_2
\end{pmatrix}=
\begin{pmatrix}
    4&7&1&0\\-6&-8&1&3
\end{pmatrix}.
$$

which provides the given set of positive Pl\"ucker coordinates. Thus, performing the inverse reordering we obtain the final result $
\begin{pmatrix}
    4&1&7&0\\-6&1&-8&3
\end{pmatrix}.
$
\end{example}

\begin{remark}
Transformations \eqref{eq:T_st} have a remarkable geometric interpretation due to the natural embedding $\mathbb{R}^n\hookrightarrow \mathbb{R}^{n+1}$ as an affine plane $x_{n+1}=1$, providing the bijection between the 
Grassmanian $G(k+1,n+1)$ and the affine Grassmaninan $Gr(k,n)$.

For instance, for the transformation $T_{3,4}$, consider the above embedding of the element $\begin{bmatrix}
   a_1&a_2&a_3&a_4\\
   b_1&b_2&b_3&b_4
\end{bmatrix}\in G(2,4)$
as the line in $\mathbb{R}^3$.

\begin{figure}[!h]
    \centering
    \includegraphics[width=0.4\textwidth]{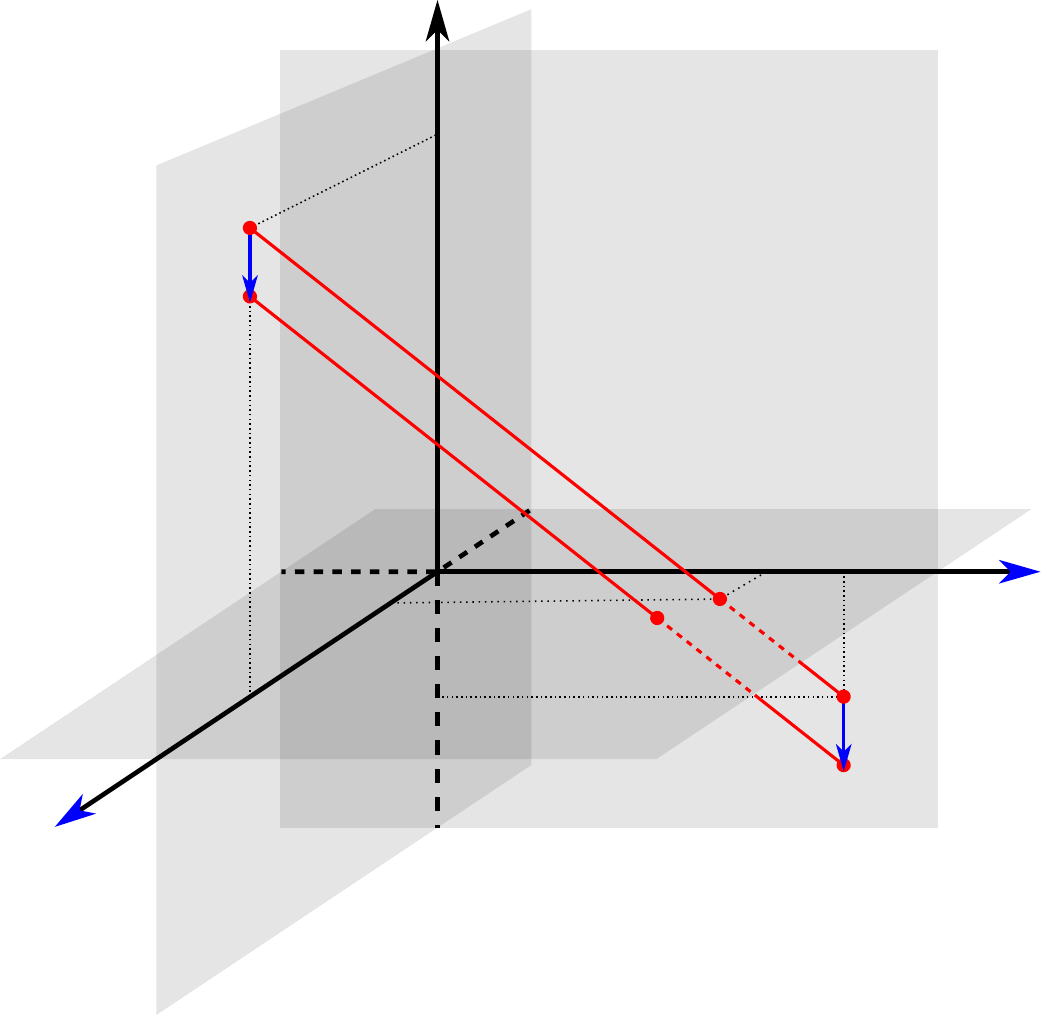}
    \caption{Affine Grassmanian}
    \label{fig:geometry}
\end{figure}

This line intersects three coordinate planes at three points: $A_1=p_{2,3}^{-1}(0,-p_{3,4},p_{2,4})$,  $A_2=p_{1,3}^{-1}(-p_{3,4},0,p_{1,4})$ and $A_3=p_{1,2}^{-1}(-p_{2,4},p_{1,4},0)$.

The transformation $T_{3,4}$ acts on these points as  
$A_1\mapsto A_1-(0,0,1)$, 
$A_2\mapsto A_2-(0,0,1)$, 
$A_3\mapsto A_3 - (-p_{23}/p_{12},p_{13}/p_{12},0)$ (see Figure \ref{fig:geometry}) In other words, transformations $T_{s,t}$ in the embedding $x_t=1$, shift the affine lines by $1$ in the direction of $e_s$.

\end{remark}

\section{Directions for further studies}
\label{Directions for further studies}

In this paper we do the first steps in the study of Euclidean subtractive algorithms for Grassmanians. 
There are numerous natural questions coming from classical Euclidean settings that are awaiting investigation. They include:

\vspace{1mm}

\begin{itemize}
\item {\bf Quality of approximation}. Classical Euclidean algorithms provide theory of best approximations for vectors in $\r^n$. What is the situation for the algorithms in Grassmanians? Is there a natural generalization of Markov numbers here? 
It would be interesting to find an analog of Littlewood conjecture.

\item {\bf Algebraic periodicity.} Which planes in $G(k,n)$ possess periodic conti\-nued-fraction sequences? IS there some analog of Lagrange theorem here?

\item{\bf Dynamic properties.} Study the analogs of Gauss maps and corresponding Gauss-Kusmin distribution. Are there some links to lattice trigonometry?

\item {\bf Generic Skew Symmetric Matrices}
 One can apply our algorithms for generic skew-symmetric matrices, not necessarily coming from the Pl\"ucker embedding of an oriented Grassmanian.

It seems that our algorithm will therefore be related to the reduction to the Hermite Normal Form skew symmetric matrices. There are numerous algorithms there, but only recently (1970th) there were some polynomial ones found, see e.g.~\cite{Grotschel1988}.

The case of Grassmanians corresponds to rank 2 matrices. 
It would be interesting to investigate the proposed algorithms for the next non-trivial case of  rank 4 matrices, coming from the fist secant variety of $G(2,n)$.

\end{itemize}

Furthermore, it would be nice to have some version of 
the following algorithms.
\begin{itemize}
\item 
An analog of Maximal Element Algorithms for $k\ge 3$.

\vspace{1mm} 
\item Lattice geometric algorithms generalising Klein polyhedra and Minkovski-Voronoi continued fractions.

\vspace{1mm}
\item Algorithms coming from algebraic number theory.
\end{itemize}

We believe that the progress in the above directions will lead to new discoveries in geometriy of Grassmanians.

\vspace{2mm}

\noindent
{\bf Acknowledgment.} The authors are grateful to K.~Khanin for useful remarks and discussions.

	\printbibliography
	
	\bigskip 
	
	\noindent
	\footnotesize \textbf{Authors' addresses:}
	
\end{document}